\documentclass[12pt,reqno]{amsart}
\usepackage{amssymb,latexsym,amsmath,amsthm,amscd,epsfig}
\usepackage[latin1]{inputenc}

\def\Im{\text{\rm Im\,}}

\def\cD{\mathcal{D}}

\def\cS{\mathcal{S}}

\def\cJ{\mathcal{J}}
\def\cP{\mathcal{P}}

\def\bpm{\begin{pmatrix}}
\def\epm{\end{pmatrix}}

\newcommand{\rf}[1]{{\eqref{#1}}}

\newcommand{\D}{{\mathbb D}}

\newcommand{\pv}{{\rm p.v.}}

\newtheorem{thm}{Theorem}[section]

\newtheorem{proposition}[thm]{Proposition}

\newtheorem{lemma}[thm]{Lemma}

\newtheorem{defn}[thm]{Definition}

\newtheorem{theorem}[thm]{Theorem}

\newtheorem{coro}[thm]{Corollary}

\theoremstyle{remark}
\newtheorem{remark}[thm]{\bf Remark}

\numberwithin{equation}{section}

\def\R{\mathbb R}
\def\C{\mathbb C}
\def\D{\mathbb D}

\def\diam{\text{diam}}

\def\dist{\operatorname{dist}}

\def\diam{\operatorname{diam}}

\newcommand{\ve}{{\varepsilon}}

\newcommand{\loc}{\operatorname{loc}}

\title{Hausdorff measure of quasicircles}

\author[I. Prause]{Istv\'an Prause}
\address{Department of Mathematics and Statistics, University of Helsinki, 
P.O.~Box 68, FIN-00014, Finland}
\email{{\tt Istvan.Prause@helsinki.fi}}

\author[X. Tolsa]{Xavier Tolsa}
\address{Instituci\'{o} Catalana de Recerca i Estudis Avan\c{c}ats (ICREA) and Departament de Matem\`{a}tiques,  Universitat Aut\`{o}noma de Barcelona, 08193 Bellaterra (Barcelona), Catalunya}
\email{{\tt xtolsa@math.uab.cat}}
\urladdr{http://mat.uab.es/~xtolsa}

\author[I. Uriarte-Tuero]{Ignacio Uriarte-Tuero}
\address{Department of Mathematics, Michigan State University, East Lansing, MI 48824, USA}
\email{{\tt ignacio@math.msu.edu}}


\thanks{I.\ P.\ was supported by projects 118634 and 1134757 of the Academy of Finland and by the Swiss NSF.
X.\ T.\ is partially supported by grants MTM2007-62817 (Spain) and 2009 SGR 420 (Catalonia). I.\ U.-T.\ 
was a postdoctoral fellow at the University of Missouri, Columbia, USA, and at Centre de Recerca Matem\`{a}tica, Barcelona, Spain, for some periods of time during the elaboration of this paper. He was partially supported by NSF grant DMS-0901524. The authors acknowledge the support of the CRM at Barcelona, where part of this paper was done 
while attending the research semester ``Harmonic Analysis, Geometric Measure Theory and Quasiconformal Mappings'' in 2009.
}

\begin{document}

\maketitle

\begin{abstract}
 S. Smirnov proved recently  \cite{Smirnov} that the Hausdorff dimension of any $K$-quasicircle is at most
$1+k^2$, where $k= (K-1)/(K+1)$. In this paper we show that if $\Gamma$ is such a quasicircle, then
$$H^{1+k^2}(B(x,r)\cap \Gamma)\leq C(k) r^{1+k^2}\qquad\mbox{for all $x\in\C$, $r>0$,}$$ 
where $H^s$ stands for the $s$-Haudorff measure. On a related note we derive a sharp weak-integrability of the derivative of the Riemann map of a quasidisk.
\end{abstract}

\section{Introduction}

A homeomorphism $f \colon \Omega \to \Omega'$ between planar domains is called \emph{$K$-quasiconformal} if it preserves orientation, belongs to the Sobolev class $W^{1,2}_{\loc}(\Omega)$ and its directional derivatives satisfy the distortion inequality
\begin{equation*}
\max_\alpha |\partial_\alpha f| \leq K \min_\alpha |\partial_\alpha f| \qquad \mbox{a.e. in $\Omega$.}
\end{equation*}
This estimate is equivalent to saying that $f$ satisfies the Beltrami equation
$$\bar \partial f(z) = \mu(z)\,\partial f(z)$$
for almost all $z\in\C$, where $\mu$ is the so called Beltrami 
coefficient or dilatation, with $\|\mu\|_\infty \leq k= (K-1)/(K+1)$.

Infinitesimally, quasiconformal mappings carry circles to ellipses with eccentricity uniformly bounded by $K$.
If $K=1$ we recover conformal maps, while for $K >1$ quasiconformal maps need not be smooth, in fact, they may distort the Hausdorff dimension of sets. The higher integrability result of Astala \cite{astalaareadistortion} provides precise estimates for this latter phenomenon. Very recently, these distortion bounds have been established even for the corresponding Hausdorff measures \cite{Lacey-Sawyer-Uriarte,tolsaqcdistortion}. In the present note we consider quasicircles and their relation to Hausdorff measures.

A $K$-{\em quasicircle} is the image of the unit circle under a $K$-quasiconformal homeomorphism of the Riemann sphere $\hat \C$. Sometimes, it will be more convenient to specialize to quasilines, these are images of the real line under a quasiconformal homeomorphism of the finite plane $\C$.
For many different characterizations of quasicircles and quasidisks (domains bounded by quasicircles), we refer the reader to \cite{gehringquasidisks}.
Complex dynamics (Julia sets, limit sets of quasi-Fuchsian groups) provide a rich source of examples of
quasicircles with Hausdorff dimension greater than one.
Astala conjectured in \cite{astalaareadistortion} that $1+k^2$ is the optimal bound on the dimension of $K$-quasicircles, where $k= (K-1)/(K+1)$. Smirnov \cite{Smirnov} proved that indeed a $K$-quasicircle has Hausdorff dimension at most $1+k^2$. 
The question of sharpness is an open problem with 
important connections to extremal behaviour of harmonic measure \cite{prause-smirnov}.
Currently, the best known lower bound appears to be the computer-aided estimate $1+0.69 k^2$
of \cite{astalarohdeschramm95}.

Our main result is the following strengthening of Smirnov's theorem in terms of $1+k^2$-dimensional Hausdorff measure $H^{1+k^2}$.
\begin{theorem}
\label{thm:main}
If $\Gamma$ is a $K$-quasicircle in $\hat \C$, then
\begin{equation*}
H^{1+k^2}(\Gamma\cap B(z,r)) \leq C(K)r^{1+k^2}\quad\mbox{ for all $z\in\C$ and with $k=\frac{K-1}{K+1}$.}
\end{equation*}
\end{theorem}

To prove this result we use a well-known factorization to ``conformal inside'' and ``conformal outside'' parts.
The conformal inside part is taken care of by Smirnov's work, we recall the necessary estimates in Section 2.
Section 3 handles the conformal outside part and we finally put together the estimates in Section 4.
Section 3 contains most of the novelties, here we adopt the technique of \cite{Lacey-Sawyer-Uriarte}
and show the boundedness of the Beurling transform with respect to some weights arising from a special packing condition. To implement the techniques in [LSUT] to our setting we have to overcome some difficulties. For instance, the arguments in [LSUT] are well suited to estimate quasiconformal distortion in terms of Hausdorff contents, while we need to obtain estimates in terms of Hausdorff measures. Also, the arguments involving  the factorization of quasiconformal maps and the use of packing conditions are much more delicate in our case.

\medskip

In Section 5 of this paper we prove another related result which
controls the expansion of the Riemann map $\phi \colon \D \to \Omega$ onto a bounded $K$-quasidisk $\Omega$. In particular, we obtain the following precise integrability condition for the derivative:
\begin{equation*}
\phi' \in \text{ weak-}L^p(\D), \quad \text{with }p=\frac{2(K^2+1)}{(K^2-1)}.
\end{equation*}
This is a strengthening of a result from \cite{prause-smirnov} to the critical exponent $p$ above and it is optimal. See Theorem \ref{thm:riemannmap} for the precise statement.

\bigskip
In the paper, as usual, the letter $C$ denotes a constant
that may change at different occurrences, while constants with subscript,
such as $C_1$ , retain their values. The notation $A \approx B$ means that there is a constant $C$ (often allowed to depend on the quasiconformality constant $K$) such that $1/C \cdot A \leq B \leq C \cdot A$. The notation $A \lesssim B$ means that there is a constant $C$ (often allowed to depend on the quasiconformality constant $K$) such that $A \leq  C \cdot B$. 
For instance, we shall frequently use this notation in conjuction with 
the well-known quasisymmetry property (see e.g.~\cite{AIM}) of quasiconformal maps. Also, as usual, if $B$ denotes a ball, $2B$ denotes the ball with the same center as $B$ and twice the radius of $B$ (and similarly for squares and other multiples.)

\section{Smirnov's theorem on quasicircles}

Before stating Smirnov's theorem, let us introduce some definitions.

\begin{defn}
 A quasiconformal mapping $f \colon \C \to \C$ is called {\em principal} if it is conformal outside some compact set $K \subset \C$ and satisfies the normalization $f(z)=z+O(1/z)$ at infinity.
\end{defn}

\begin{defn}
 A quasiconformal mapping $f \colon \C \to \C$ is called {\em antisymmetric} (with respect to the real line) (or equivalently we say that it has antisymmetric dilatation $\mu$)
if its Beltrami coefficient $\mu$ satisfies
\[
 \mu(z)=-\overline{\mu(\bar z)} \qquad \text{for a.e. $z \in \C$}.
\]
\end{defn}
The significance of this definition is that by making use of a symmetrization procedure \cite{Smirnov} shows that, any $K$-quasiline may be represented as the image of the real line under a $K$-quasiconformal antisymmetric map.

Now we state Smirnov's bound on the dimension of quasicircles.
\begin{theorem}[\cite{Smirnov}]
The Hausdorff dimension of a $\frac{1+k}{1-k}$-quasicircle is at most $1+k^2$, for any $k\in[0,1)$.
\end{theorem}

We shall need the following formulation. This is implicit in \cite{Smirnov} and 
appears exactly as stated in \cite[Theorem 13.3.6]{AIM}.

\begin{theorem}\label{teosmirnov}
Let $f:\C\to\C$ be a $\frac{1+k}{1-k}$-quasiconformal map, with $0<k<1$. Suppose
that $f$ has antisymmetric dilatation and that it is principal and 
conformal outside the unit disk.
Let $B_j=B(z_j,r_j)$, $1\leq j\leq n$, be a collection of pairwise disjoint disks contained in the unit disk such that $z_j\in\R$,
so that $f$ is conformal on each $B_j$.
For $0<t\leq 2$, let $t(k)$ be such that
$$\frac1{t(k)}- \frac12 = \frac{1-k^2}{1+k^2}\,\left(\frac1t - \frac12\right).$$
Then,
$$\biggl(\sum_{j=1}^n (|f'(z_j)|\,r_j)^{t(k)} \biggr)^{\frac1{t(k)}}
\leq 8 \biggl( \sum_{j=1}^n r_j^t 
\biggr)^{\frac{1-k^2}{1+k^2} \frac1t}.$$
\end{theorem}

As an easy corollary we give the special case $t=1$ in an invariant form.

\begin{coro}\label{corosmirnov}
Let $f:\C\to\C$ be a $\frac{1+k}{1-k}$-quasiconformal map, with $0<k<1$, with
antisymmetric dilatation.
Let $B_j$, $1\leq j\leq n$, be a collection of pairwise disjoint disks 
contained in a disk $B$.
Suppose that the disks $B_j$, $B$,  are centered on the real line and that $f$ is conformal on the disks $B_j$.
Then,
$$\frac{\sum_{j=1}^n \diam(f(B_j))^{1+k^2}}{\diam(f(B))^{1+k^2}} 
\leq C(k) \Biggl( \frac{\sum_{j=1}^n \diam(B_j)}{\diam(B)}
\Biggr)^{1-k^2}.$$
In particular, notice that in the above situation,
$$\sum_{j=1}^n \diam(f(B_j))^{1+k^2} \leq C(k) \diam(f(B))^{1+k^2}.$$
\end{coro}

\begin{proof}
We factorize $f=f_2\circ f_1$, where $f_1$, $f_2$ are both $K$-quasiconformal
maps, with $f_1$ principal and conformal on $\C\setminus 2B$, and
$f_2$ is conformal on $f_1(2B)$. 
If we denote by $\mu_f$ the dilatation of $f$, we assume that the one of $f_1$ is
$\mu_{f_1}= \chi_{2B}\mu_f$.
Let $g(z)=az+ b$ be the function that maps the unit disk to $2B$. 
The function $h=g^{-1}\circ f_1\circ g$ with the collection of disks $\{g^{-1}(B_j) \}$ verifies the assumptions of Theorem \ref{teosmirnov}, and so specializing to $t=1$ we obtain
$$\sum_{j=1}^n \diam(h(g^{-1}(B_j)))^{1+k^2} 
\leq C(k) \Bigl(\sum_{j=1}^n \diam(g^{-1}(B_j))\Bigr)^{1-k^2}.$$
Since $\diam(h(g^{-1}(B_j)))= \diam(g^{-1}(f_1 (B_j))) = a^{-1}\diam(f_1 (B_j))$, we deduce
$$\sum_{j=1}^n \diam(f_1 (B_j))^{1+k^2} 
\leq C(k) a^{2k^2} \Bigl(\sum_{j=1}^n \diam(B_j)\Bigr)^{1-k^2}.$$
On the other hand, since $f_2$ is conformal on $f_1(2B)$, by 
Koebe's distortion theorem and quasisymmetry, 
$$\frac{\diam(f_2(f_1(B_j)))}{\diam(f_2(f_1(2B)))} \approx 
\frac{\diam(f_1(B_j))}{\diam(f_1(2B))}.$$
Therefore,
$$\sum_{j=1}^n \diam(f(B_j))^{1+k^2} 
\leq C(k) \Bigl(\sum_{j=1}^n \diam(B_j)\Bigr)^{1-k^2}\,
\left(\frac{\diam(f(2B))}{\diam(f_1(2B))}\right)^{1+k^2}\,
a^{2k^2}.$$
Replacing $a=\diam(B)$ and taking into account that $\diam(f_1(2B))\approx\diam(2B)$
by Koebe's distortion theorem, the corollary follows.
\end{proof}

\section{Estimates for the ``conformal outside'' map}

\subsection{Smooth and packed families of dyadic squares}\label{secsmooth}

Let us denote the family of dyadic squares by $\cD$. Let $0<t<2$ be fixed. 
Let $\{Q\}_{Q\in\cJ}$ a family of pairwise disjoint dyadic squares. Given $\tau\geq1$, 
we say that $\cJ$ is a {\em $\tau$-smooth}
family if 
\begin{itemize}
\item[(a)] If $P,Q\in\cJ$ are such that $2P\cap 2Q\neq\varnothing$, then
$\tau^{-1}\ell(P)\leq \ell(Q)\leq \tau\ell(P).$

\item[(b)] $\sum_{Q\in\cJ} \chi_{2Q}\leq \tau.$
\end{itemize}
We say that $\cJ$ is $\alpha$-packed if for any dyadic square $R\in\cD$ which contains at least two
squares from $\cJ$,
\begin{equation}\label{eqbb}
\sum_{Q\in\cJ:\,Q\subset R} \ell(Q)^t\leq \alpha\,\ell(R)^t.
\end{equation}
At first sight, the fact that we ask \rf{eqbb} only for squares $R$ which contain at least two squares
from $\cJ$ may seem strange. In fact, \rf{eqbb} holds with $\alpha=1$ for any square $R$ which contains
a unique square $Q$. The advantage of the formulation above is that for $\alpha$ arbitrarily small, there exist
$\alpha$-packed families $\cJ$, which is not the case if we allow $R$ to contain a unique square $Q\in\cJ$.
This fact  will play an important role below.

Recall that a quasisquare is the image of a square by a 
quasiconformal map. 
Given a quasiconformal map $f$, we say that $Q$ is a (dyadic) $f$-quasisquare
if it is the image of a (dyadic) square by $f$. 
If $Q=f(P)$ is an $f$-quasisquare, we denote $aQ:=f(aP)$.
If now $\{Q\}_{Q\in\cJ}$ is family of pairwise disjoint dyadic $f$-quasisquares, we say
that it is $\tau$-smooth and $\alpha$-packed if it verifies the properties above, replacing
$\ell(P),\ell(Q),\ell(R)$ by $\diam(Q),\diam(P),\diam(R)$, and \rf{eqbb} is required to hold for all
dyadic $f$-quasisquares $R$ containing at least two quasisquares from $\cJ$.


\subsection{Boundedness of the Beurling transform with respect to some weights}

In the following, given a non negative measurable function (i.e. a weight) $\omega$
and a subset $A\subset\C$, we use the standard notation
$$\omega(A) = \int_A \omega\,dm.$$
Recall also that the Beurling transform of a function $f:\C\to\C$ is given by
$$\cS f(z) = \frac{-1}{\pi}\,\pv \int_\C \frac{f(\xi)}{(z-\xi)^2}\,dm(\xi).$$

\begin{proposition}\label{propo9}
Let $0<t<2$.
Let $\cP=\{ P_i\}_{i=1}^N$ be a finite $\tau$-smooth $\alpha$-packed 
family of pairwise disjoint dyadic squares.
Denote $\overline P = \bigcup_{i=1}^N P_i$ and set
$$\omega = \sum_{P\in\cP} \ell(P)^{t-2}\,\chi_P.$$
Then, the Beurling transform is bounded in $L^2(\omega)$. That is to say,
\begin{equation}\label{eqacot}
\|\cS(f\chi_{\overline P})\|_{L^2(\omega)}\leq C_1\|f\|_{L^2(\omega)}
\quad \mbox{for all $f\in L^2(\omega)$}.
\end{equation}
The constant $C_1$ only depends on $t$, $\tau$, and $\alpha$.
\end{proposition}

\begin{proof}
We will follow the arguments of \cite{Lacey-Sawyer-Uriarte} quite closely. By interpolation, it is  enough to show that
$$\int_G |\cS(\chi_F)|\,\omega\,dm\leq C_p\,\omega(F)^{1/p}\,\omega(G)^{1/p'}
\quad \mbox{ for all $F,G\subset\overline P$ and $1<p<\infty$}.$$
To prove this estimate, for $f\in L^1_{loc}(\C)$, we split $\chi_{\overline P}\cS(f\chi_{\overline P})$ into a local and a non local part
as follows
$$\chi_{\overline P}\cS(f\chi_{\overline P}) = \sum_{i=1}^N\chi_{\overline P\cap 2 P_i}\cS(f\chi_{P_i}) + 
\sum_{i=1}^N\chi_{\overline P\setminus 2 P_i}\cS(f\chi_{P_i}) = :\cS_{local}(f) +
\cS_{non}(f).$$
For the local part we will use the boundedness of $\cS$ in $L^p(\C)$, the fact that $\omega$ is a constant times Lebesgue measure on each $P_i$, and moreover,
$\omega\approx\omega_{|P_i}$ on $2 P_i \cap \overline P$, because of the property (a) in Subsection \ref{secsmooth}. Then,
 by H\"older's inequality,
\begin{align*}
\int_G |\cS_{local}(\chi_F)|\,\omega\,dm & \leq \sum_i 
\int_{2 P_i\cap G} |\cS(\chi_{P_i\cap F})|\,w\,dm\\
& \leq C_p\sum_i|P_i\cap F|^{1/p}\,|2P_i\cap G|^{1/p'}\,\omega_{|P_i}\\
& \approx  C_p\sum_i\omega(P_i\cap F)^{1/p}\,\omega(2 P_i\cap G)^{1/p'}\\
& \leq C_p\Bigl(\sum_i\omega(P_i\cap F)\Bigr)^{1/p}\,\Bigl(\sum_i\omega(2 P_i\cap G)\Bigr)^{1/p'} \\ &\leq C\,\omega(F)^{1/p}\omega(G)^{1/p'},
\end{align*}
where we used the property (b) of smooth families in the last inequality.

Consider now the non local part. 
For any two squares $P,Q$, denote
$$D(P,Q)=\dist(P,Q) + \ell(P) + \ell(Q).$$
Since for all $P,Q\in\cP$ such that $Q\setminus 2 P\neq\varnothing$,
$$\dist(P,Q\setminus 2 P)\approx D(P,Q),$$
for $x\in P\in\cP$ we have
$$|\cS_{non}f(x)|\lesssim \sum_{Q\in\cP:Q\neq P}\frac1{D(P,Q)^2}\int_Q |f|\,dm =:Tf(x).$$ 
We will show that
\begin{equation}\label{eqtt}
\int_G |T(\chi_F)|\,\omega\,dm\leq C_p\,\omega(F)^{1/p}\,\omega(G)^{1/p'}.
\end{equation}
Let $M_\omega$ be the following maximal operator:
$$M_\omega f(x) =\sup_{x\in Q}\frac1{\ell(Q)^t} \int_Q |f|\,\omega\,dm,$$
where the supremum is taken over all the squares containing $x$. Since
$\omega(Q)\leq C \ell(Q)^t$ for any square $Q$ (by the packing condition), 
it follows by standard arguments (using covering lemmas) that 
$M_\omega$
is bounded in $L^p(\omega)$, $1<p\leq\infty$, and from $L^1(\omega)$ 
to $L^{1,\infty}(\omega)$. To prove \rf{eqtt}, for a fixed $C_1>2$ we define
$$F' =\left\{ \begin{array}{l} F\quad \mbox{if $8\omega(F)\leq\omega(G)$,}\\
\mbox{}\\
F\cap \{ M_\omega(\chi_G)\leq C_1\omega(G)/\omega(F)\}\quad\mbox{otherwise.}
\end{array}
\right.$$
From the weak $(1,1)$ inequality for $M_\omega$, we have $\omega(F')\geq \omega(F)/2$ 
if $C_1$ is chosen big enough.
As in \cite{Lacey-Sawyer-Uriarte}, to prove \rf{eqtt} it is enough to show that
$$\int_G (T\chi_{F'})\,\omega\,dm\leq C\,\min\bigl(\omega(F),\,\omega(G)\bigr),$$
(iterating this estimate, \rf{eqtt} follows). 
We have
\begin{align*}
\int_G (T\chi_{F'})\,\omega\,dm & = 
\sum_{P\in\cP}\sum_{Q\in\cP:P\neq Q}\frac{|F'\cap Q|}{D(P,Q)^2}\,
\omega(G\cap P) \\ &= \sum_{Q\in\cP}|F'\cap Q|\sum_{P\in\cP:P\neq Q}\frac{1}{D(P,Q)^2}\,
\omega(G\cap P).
\end{align*}
Denote by $A_k(Q)$ the family of squares $P\in\cP$ such that 
$$2^k\ell(Q)\leq D(P,Q)<2^{k+1}\ell(Q).$$
Then,
\begin{align*}
\sum_{P\in\cP:P\neq Q}\frac{1}{D(P,Q)^2}\,
\omega(G\cap P) & \leq C \sum_{k\geq0}\frac1{\ell(2^kQ)^2}\sum_{P\in A_k(Q)}
\omega(G\cap P)\\ & \leq C\sum_{k\geq0}\frac1{\ell(2^kQ)^2} \, \omega(G\cap 2^{k+3}Q),
\end{align*}
because the squares $P\in A_k(Q)$ are contained $2^{k+3}Q$.
By the packing condition \rf{eqbb}, we have
$$\omega(G\cap 2^{k+3}Q)\leq C\,\ell(2^{k+3}Q)^t\leq C\,\ell(2^kQ)^t.$$
On the other hand, assuming $Q\cap F'\neq\varnothing$, by the definition of $F'$,
$$\omega(G\cap  2^{k+3}Q) \leq C\ell(2^{k+3}Q)^t\,\frac{\omega(G)}{\omega(F)}\leq C\,\ell(2^kQ)^t\,
\frac{\omega(G)}{\omega(F)}.$$
So we get,
$$\omega(G\cap 2^{k+3}Q) \leq C\,\min\biggl(1,\frac{\omega(G)}{\omega(F)}\biggr)\,\ell(2^kQ)^t
=:C A\,\ell(2^kQ)^t.$$
Therefore, since $t-2<0$,
\begin{align*}
\int_G (T\chi_{F'})\,\omega\,dm &\leq C A\sum_{Q\in\cP}|F'\cap Q|
\sum_{k\geq0}\ell(2^kQ)^{t-2}\leq C A\sum_{Q\in\cP}|F'\cap Q|\,
\ell(Q)^{t-2}\\ 
&= CA\omega(F')\leq C\min\bigl(\omega(F),\,\omega(G)\bigr).
\end{align*}
\end{proof}

\begin{remark}
Given a $K$-quasiconformal map $f:\C\to\C$,
the proposition above also holds (by the same proof) if we consider a family of dyadic $f$-quasisquares instead of dyadic
squares, with $$\omega=\sum_{P\in\cP} \diam(P)^{t-2}\,\chi_P.$$
 Then, the constant $C_1$ in \rf{eqacot} depends on $K$, besides $\alpha$, $\tau$ and $t$.
\end{remark}

\bigskip



\subsection{Distortion when $f$ is conformal outside a family of quasisquares}

\begin{lemma} \label{lemkpetit}
Let $0<t<2$.
Let $\{Q\}_{Q\in\cJ}$ be a finite $\tau$-smooth $\alpha$-packed family of pairwise disjoint dyadic 
$g$-quasisquares, where $g$ is some $K_0$-quasiconformal map with $K_0\leq M_0$, $\alpha\leq 1$. 
Denote $F=\bigcup_{Q\in\cJ} Q$ and let  $f:\C\rightarrow\C$ be a principal $K$-quasiconformal map
 conformal on $\C\setminus F$.
There exists $\delta_0 =\delta_0(t,\tau,M_0)>1$ such that if $K\leq \delta_0$, then
 $$\sum_{Q\in\cJ}\diam(f(Q))^t\leq C_2\sum_{Q\in\cJ}\diam(Q)^t,$$
 with $C_2$ depending on only on $t,\tau,M_0$. 
\end{lemma}

The proof is analogous to the one in \cite[Lemma 5.6]{Lacey-Sawyer-Uriarte}, using the fact the the Beurling transform is bounded in $L^2(\omega)$, by Proposition \ref{propo9}.

We wish to extend the preceding result to $K$-quasiconformal maps with $K$ arbitrarily large. As usual, 
we will do this by an appropriate factorization of $f$. First we need the following technical result.

\begin{lemma} \label{lemkgran1}
Let $0<t<2$.
Let $\{Q\}_{Q\in\cJ}$ be a finite $\tau$-smooth $\alpha$-packed family of pairwise disjoint dyadic 
$g$-quasisquares, where $g$ is some $K_0$-quasiconformal map with $K_0\leq M_0$, $\alpha\leq 1$. 
Denote $F=\bigcup_{Q\in\cJ} Q$ and let  $f:\C\rightarrow\C$ be a principal $K$-quasiconformal map
 conformal on $\C\setminus F$.
There exists $\delta_0=\delta_0(t,\tau,M_0)>1$ such that if $K\leq \delta_0$, then
for any $g$-dyadic quasisquare $R$,
 $$\sum_{Q\in\cJ:Q\subset R}\diam(f(Q))^t\leq C_3\,\frac{\sum_{Q\in\cJ:Q\subset 3R}\diam(Q)^t}{\diam(R)^t}\,
 \diam(f(R))^t ,$$
 with $C_3=C_3(t,\tau,M_0)$. In particular, 
 the family of $(f\circ g)$-quasisquares 
$\{f(Q)\}_{Q\in\cJ}$ is $(9C_3\alpha)$-packed.
\end{lemma}

\begin{proof}
We factorize $f=f_2\circ f_1$, where $f_1$, $f_2$ are $K$-quasiconformal
maps, with $f_1$ conformal in $\bigl(\C\setminus \bigcup_{Q\in\cJ} Q\bigr)\cup (\C\setminus 3R)$, and
$f_2$ is conformal on $f_1(3R)$. By Lemma \ref{lemkpetit}, we have
\begin{equation}\label{eq41}
\sum_{Q\in\cJ:Q\subset 3R}\diam(f_1(Q))^t\leq 
C_2\sum_{Q\in\cJ:Q\subset 3R}\diam(Q)^t.
\end{equation}
By Koebe's distorsion theorem and quasisymmetry, 
since $f_2$ is conformal in $f_1(3R)$,
for every $Q\subset R$,
$$\frac{\diam(f_2(f_1(Q)))}{\diam(f_2(f_1(3R)))} \approx
\frac{\diam(f_1(Q))}{\diam(f_1(3R))}.$$
Thus,
$$\sum_{Q\in\cJ:Q\subset R}\frac{\diam(f(Q))^t}{\diam(f(3R))^t}
\approx \sum_{Q\in\cJ:Q\subset R}\frac{\diam(f_1(Q))^t}{\diam(f_1(3R))^t}.$$
The lemma follows from this estimate, \rf{eq41}, and the fact that
$\diam(f_1(3R))\approx\diam(3R)$, since $f_1$ is principal and conformal on 
$\C\setminus 3R$.
\end{proof}

\begin{lemma} \label{lemkgran2}
Let $0<t<2$.
Let $\{Q\}_{Q\in\cJ}$ be a finite $\tau$-smooth $\alpha$-packed family of pairwise disjoint dyadic 
$g$-quasisquares, where $g$ is some $K_0$-quasiconformal map with $K_0\leq M_0$. 
Denote $F=\bigcup_{Q\in\cJ} Q$ and let  $f:\C\rightarrow\C$ be a principal $K$-quasiconformal map
 conformal on $\C\setminus F$ with $K\leq M_0$.
There exists $\delta_1>0$ small enough depending only on $t,\tau,M_0$ such that if $\alpha\leq \delta_1$, then
 $$\sum_{Q\in\cJ}\diam(f(Q))^t\leq C_4\sum_{Q\in\cJ}\diam(Q)^t,$$
 with $C_4$ depending on only on $t,\alpha,\tau,M_0$. 
\end{lemma}

\begin{proof} 
Notice that for any $K'$-quasiconformal map $h$ with $K'\leq M_0$, the family $\{h(Q)\}_{Q\in\cJ}$ is 
$\tau'$-smooth, with $\tau'$ depending only on $\tau$ and $M_0$. Let $n$ be big enough
so that $M_0^{1/n}\leq \delta_0'$, where $\delta_0'=\delta_0(t,\tau',M_0)$ (with $\delta_0$ from Lemma \ref{lemkpetit}).

We factorize $f=f_n\circ\cdots \circ f_1$
so that each $f_i$ is $K^{1/n}$-quasiconformal on $\C$, and moreover  $f_i$ is conformal in 
$\C\setminus \bigcup_{Q\in\cJ} f_{i-1}\circ \cdots \circ f_1(Q)$, for $i\geq 1$ (with $f_0=id$). Notice that for all $i$, the quasisquares
$f_i\circ \cdots \circ f_1(Q)$ are $\tau'$-smooth (since
$f_i\circ \cdots \circ f_1$ is $K_i'$-quasiconformal with $K_i'\leq M_0$).
Suppose that $\alpha$ is small enough so that
$$(C_3')^{n}\alpha\leq 1,$$
where $C_3' = 9C_3(t,\tau',M_0)$.
Then, Lemma \ref{lemkgran1} can be applied repeatedly to deduce 
that for each $i\leq n$
the family of quasisquares $\{f_i\circ \cdots \circ f_1(Q)\}_{Q\in\cJ}$ is $(C_3')^i\alpha$-packed, and thus
$1$-packed (without loss of generality, we assume $C_3'\geq1$). 
As a consequence, Lemma \ref{lemkpetit} can also be applied repeatedly 
to get
$$\sum_{Q\in\cJ}\diam(f(Q))^t\leq C_2\sum_{Q\in\cJ}\diam(f_{n-1}\circ \cdots \circ f_1(Q))^t \leq \cdots 
\leq C_2^n\sum_{Q\in\cJ}\diam(Q)^t.$$
\end{proof}


\section{Gluing ``conformal inside'' and ``conformal outside''}

In the following lemma we make use of the conformal inside vs.~outside decomposition.

\begin{lemma}\label{mainlem}
Let $f:\C\to\C$ be a $\frac{1+k}{1-k}$-quasiconformal map with antisymmetric dilatation.
 Let $\{Q\}_{Q\in\cJ}$ be a finite family of pairwise disjoint squares with 
equal side length
centered on $\R$ which are contained in another square $R$ centered on $\R$. Then, 
\begin{equation}\label{eqindu1}
\sum_{Q\in\cJ} \diam(f(Q))^{1+k^2} \leq 
C(k)\left(\frac{\sum_{Q\in\cJ} \ell(Q)}{\ell(R)}\right)^{1-k^2}
\diam(f(R))^{1+k^2}.
\end{equation}
\end{lemma}

\begin{proof}
First, let us make the assumption that $f$ is principal. We are going to relax this assumption at the end of the proof.
Using quasisymmetry if necessary, 
we may assume that the squares $Q\in\cJ$ 
belong to a dyadic lattice (a translation of the usual dyadic lattice, say).
Take a small constant $0<\alpha<1$ to be fixed below.
We will prove \rf{eqindu1} assuming that $\cJ$ is $\alpha$-packed.
The general statement follows easily from this particular case: since the squares $Q\in\cJ$ have equal side 
length and are centered on $\R$, we can easily split $\cJ=\cJ_1\cup\ldots \cup \cJ_m$ so that $m\approx
1/\alpha$ and each $\cJ_i$ is $\alpha$-packed (recall that the constant $\alpha$ in the definition of 
$\alpha$-packings only involves squares $R$ which contain at least two squares $Q\in\cJ$). Then we apply 
\rf{eqindu1} to each family $\cJ_i$ and we add the resulting estimates. This introduces an additional
multiplicative constant $m$ on the right hand side of  \rf{eqindu1}. We note here that
$m\approx 1/\alpha$ and we will later choose $\alpha$ depending only on $K=\frac{1+k}{1-k}$.

So assume that the squares $Q\in\cJ$ are $\alpha$-packed and denote $V=\bigcup_i Q_i$.
Take a decomposition
$f=f_2\circ f_1$, where $f_1$, $f_2$ are principal $K$-quasiconformal
mappings. We require $f_1$ to be conformal in $V$ and $f_2$
outside $f_1(\overline{V})$. Further, we suppose that the 
dilatation of $f_1$ is $\mu_{f_1}=\chi_{\C\setminus V}\mu_f$
so that $\mu_{f_1}$ is also antisymmetric.

By Corollary \ref{corosmirnov}, for any square $P$ centered on $\R$,
\begin{equation}\label{eqpack4}
\sum_{Q\in\cJ:Q\subset  P}\diam(f_1(Q))^{1+k^2}\leq C(k)
\left(\frac{\sum_{Q\in\cJ:Q\subset P} \ell(Q)}{\ell(P)}\right)^{1-k^2}\diam(f_1(P))^{1+k^2}.
\end{equation}
In particular, the family of quasisquares $\{f_1(Q)\}_{Q\in\cJ}$ is $C(k) \alpha^{1-k^2}$-packed.
It is also clear that they form a $\tau$-smooth family, with $\tau$ depending
only on $K$.
Therefore, if $\alpha$ has been chosen small enough (depending only on $K$), from Lemma \ref{lemkgran2} 
we deduce that
$$\sum_{Q\in\cJ}\diam(f_2(f_1(Q)))^{1+k^2}\leq C\sum_{Q\in\cJ}\diam(f_1(Q))^{1+k^2},$$
and so by \rf{eqpack4} with $P=R$,
$$\sum_{Q\in\cJ}\diam(f(Q))^{1+k^2}\leq C(k)
\left(\frac{\sum_{Q\in\cJ} \ell(Q)}{\ell(R)}\right)^{1-k^2}\diam(f_1(R))^{1+k^2}.$$

Since $f_2$ is principal and conformal outside $f_1(R)$, by Koebe's distortion theorem we deduce
$\diam(f(R))\approx \diam(f_1(R))$, and
thus \rf{eqindu1} follows for a principal mapping.
We reduce the general case to this one. Suppose that $f$ is not necessarily principal antisymmetric map.
Take a decomposition $f=g_2 \circ g_1$, where $g_1$ is principal antisymmetric $K$-quasiconformal map which is conformal outside $3R$ and $g_2$ is a $K$-quasiconformal map which is conformal on $g_1(3R)$.
This decomposition is analogous to the one used in Corollary \ref{corosmirnov} and again by Koebe's distortion theorem and quasisymmetry we have for every $Q\in\cJ$
\[ \frac{\diam(f(Q))}{\diam(f(3R))} \approx \frac{\diam(g_1(Q))}{\diam(g_1(3R))}.
\]
Now the lemma follows from \rf{eqindu1} applied to the principal map $g_1$.
\end{proof}

\begin{theorem}\label{mainteo}
Let $f:\C\to\C$ be an antisymmetric $\frac{1+k}{1-k}$-quasiconformal map. Then, for any compact subset
$E\subset \R$ and any ball $B\subset\C$ centered on $\R$ which contains $E$,
$$H^{1+k^2}(f(E))\leq C(k)\diam(f(B))^{1+k^2}\,\left(\frac{H^1(E)}{\diam(B)}\right)^{1-k^2}.$$
\end{theorem}

\begin{proof}
Consider an arbitrary covering $E\subset\bigcup_i I_i$ by a finite number of 
pairwise different
dyadic intervals of length $\ve\diam(B)$. 
Consider squares concentric with $I_i$ with $\ell(Q_i)=\diam(I_i)$. By Lemma
\ref{mainlem} and quasisymmetry we  deduce that
$$
\sum_{i} \diam(f(I_i))^{1+k^2} \leq 
C(k)\left(\frac{\sum_{i} \diam(I_i)}{\ell(B)}\right)^{1-k^2}
\diam(f(B))^{1+k^2}.$$
Because of the H\"older continuity of quasiconformal maps (see e.g.~\cite{AIM}), for each $i$, with a constant $C_5$ depending on $k$,
$$\frac{\diam(f(I_i))}{\diam(f(B))}\leq C_5 \left(
\frac{\diam(I_i)}{\diam(B)}\right)^{\frac{1-k}{1+k}} = C_5\ve^{\frac{1-k}{1+k}}.$$
Therefore, with $\delta=C_5\ve^{\frac{1-k}{1+k}}$
$$H^{1+k^2}_{\delta} (f(E)) \leq 
C(k)\left(\frac{\sum_{i} \diam(I_i)}{\diam(B)}\right)^{1-k^2}
\diam(f(B))^{1+k^2}.$$
By the definition of length on $\R$, we have $\sum_{i} \diam(I_i)\leq H^1(U_\ve(E)\cap \R)$,
where $U_\ve$ stands for the $\ve$-neighborhood, and we assume that $I_i\cap E\neq\varnothing$. Thus,
$$H^{1+k^2}_{\delta}(f(E)) \leq 
C(k)\left(\frac{H^1(U_{\ve\diam(B)}(E)\cap \R)}{\diam(B)}\right)^{1-k^2}
\diam(f(B))^{1+k^2}.$$
Letting $\ve\to0$, the theorem follows.
\end{proof}

Now we are ready to prove Theorem \ref{thm:main} stated in the introduction.

\begin{proof}[Proof of Theorem \ref{thm:main}]
If $\Gamma$ is a $K$-quasiline, then $\Gamma = f(\R)$ with some $K$-quasiconformal map $f \colon \C \to \C$.
As we remarked in Section 2, we may further suppose that $f$ is antisymmetric. Our goal is to show
\begin{equation}
\label{eq:main}
H^{1+k^2}(\Gamma\cap B(z,r)) \leq C(k)r^{1+k^2}\quad\mbox{ for all $z\in\C$.}
\end{equation}
First of all, we may assume that $z \in \Gamma$. Indeed, if $\Gamma \cap B(z,r)=\emptyset$ then there is nothing to prove. Otherwise, we can find $z_0 \in \Gamma$
such that $B(z,r) \subset B(z_0,2r)$ and hence by replacing $B(z,r)$ with a twice larger ball
we may assume that the center lies on $\Gamma$.
Let us set $E=f^{-1}(\bar B \cap \Gamma) \subset \R$. Using quasisymmetry we can easily find a ball $B_0$ centered on $\R$ which contains $E$ and such that $\diam f(B_0) \approx r$.
We apply now Theorem \ref{mainteo} with $E$ and $B_0$ as above and find that \eqref{eq:main} holds true.

The case where $\Gamma$ is a quasicircle in $\C$ can be reduced to the one of a quasiline. Indeed, with the help of a M\"obius transformation we may pass to a quasiline and see that \eqref{eq:main} holds, for instance, with $r \leq \diam \Gamma /10$. For $r > \diam \Gamma/10$ we just use \eqref{eq:main} for a finite number of balls of radius $\diam \Gamma /10$.
\end{proof}


\section{Boundary expansion of the Riemann map}

The Riemann map $\phi \colon \D \to \Omega$ onto a quasidisk is H\"older continuous up to the boundary, in short, quasidisks are H\"older domains.
For a map with $K$-quasiconformal extension the sharp H\"older exponent is $1-k$ \cite{pommerenke75}, where $k=(K-1)/(K+1)$, as usual.
Very recently, the following counterpart was established in terms of the integrability of the derivative.

\begin{theorem}[{\cite[Corollary 3.9]{prause-smirnov}}]
If $\phi \colon \D \to \C$ is a conformal map with $K$-quasiconformal extension then
\[ \phi' \in L^p(\D) \quad \mbox{for all } \quad 2 \leq p < \frac{2(K+1)}{K-1}.
\]
The upper bound for the exponent is the best possible.
\end{theorem}

In next theorem we prove the 
 weak-integrability of $\phi'$ in the borderline case
 $p=\frac2 k = \frac{2(K+1)}{K-1}$.
In terms of area distortion for subsets of the unit disk, the exponent $1/K$ from Astala's theorem improves to $1-k$, just as
the H\"older continuity exponent does.

\begin{theorem}
\label{thm:riemannmap}
If $\phi \colon \D \to \C$ is a conformal map with $K$-quasiconformal extension to $\C$ then
$\phi' \in \text{weak-}L^{2/k}(\D)$ with $k=(K-1)/(K+1)$. More precisely,
\begin{equation} 
 \label{eq:weakintegrability}
\left| \{ z \in \D : |\phi'(z)| > \rho \} \right| \leq C(K) |\phi'(0)|^{2/k} \rho^{-2/k}
\quad \text{for any $\rho >0$}.
\end{equation} 
In terms of area distortion, for any Borel set $E \subset \D$, 
\[ |\phi(E)| \leq C(K) |\phi'(0)|^2 |E|^{1-k}.
\]
\end{theorem}

\begin{remark}
The power map $z \mapsto z^{1-k}$ maps conformally the upper half-plane to a sector domain of angle $(1-k)\pi$
and admits a $\frac{1+k}{1-k}$-quasiconformal extension to $\C$ \cite{becker-pommerenke98}. 
This example shows that Theorem \ref{thm:riemannmap} is sharp up to the numerical value of the constant
terms involved.
\end{remark}

First we will prove the following lemma, as an application of 
the theorems of Smirnov and Astala.
 
\begin{lemma}
\label{lem:areaexpansion}
Let $\psi: \C \to \C$ be a principal $K$-quasiconformal map which is conformal in $\C_+=\{z: \Im z >0 \}$ and outside $\D$. 
Let $B_j=B(z_j,r_j)$, $1\leq j\leq n$, be a collection of pairwise disjoint disks contained in the unit disk such that $z_j\in\R$. We set $E=\cup B_j$.
Then we have the following estimate for area expansion
\begin{equation}
 \label{eq:areaexpansion}
|\psi(E)| \leq C(K) |E|^{1-k},
\end{equation}
with $k=(K-1)/(K+1)$.
\end{lemma}

\begin{proof}
In order to deduce \eqref{eq:areaexpansion} it is sufficient to assume that
$\psi$ is conformal in $E$. Otherwise, we use the usual decomposition to principal $K$-quasiconformal mappings $\psi=\psi_2 \circ \psi_1$ where $\psi_1$ is conformal on $E$ and $\psi_2$ conformal outside $\psi_1(E)$. Now invoking the fact $|\psi_2(\psi_1(E))| \leq K |\psi_1(E)|$ from \cite{astalaareadistortion} matters are reduced to $\psi_1$. In the rest of the proof we assume that
$\psi$ is conformal in $E$ and for notational convenience replace $K$ by $K^2$, that is, 
$\psi$ assumed to be globally $K^2$-quasiconformal. The exponent $1-k$ in \eqref{eq:areaexpansion} then takes the form
\begin{equation}
\label{eq:1-k} 
 1-\frac{K^2-1}{K^2+1} = \frac{2}{K^2+1}=\frac1K \cdot \frac{1-k^2}{1+k^2},
\end{equation}
where $k=(K-1)/(K+1)$.

We use the symmetrization result of \cite{Smirnov}: $\psi$ can be written as a superposition $\psi=f \circ g$ of a $K$-quasiconformal map $g$ symmetric with respect to $\R$ followed by a $K$-quasiconformal antisymmetric map $f$. Both of these maps are normalized to be principal. 
We observe that under the conformality assumptions on $\psi$, the map $g$ is conformal in the disks $B_j$ and outside $\D$ while $f$ is conformal in the quasidisks $g(B_j)$ and outside $g(\D)$.
In view of Koebe's $1/4$ theorem,
\[ \hat{B_j}:=B\left( g(z_j),\frac14 |g'(z_j)| r_j \right) \subset g(B_j) \quad \text{and} \quad
g(\D) \subset B(0,2).
\]
We apply Theorem \ref{teosmirnov} with $t=2$ for the map $f$ and disks $\hat{B_j} \subset B(0,2)$,
\[
 \sum_{j=1}^n \left( |f'(g(z_j))| \, |g'(z_j)| r_j \right)^2 \leq 
C \left( \sum_{j=1}^n (|g'(z_j)| r_j)^2 \right)^{\frac{1-k^2}{1+k^2}}.
\]
As $\psi'(z_j)=f'(g(z_j))g'(z_j)$, we may write the previous inequality as the comparison of area
\[ 
 |\psi(E)| \approx \sum_{j=1}^n (|\psi'(z_j)| r_j )^2 \leq C(K) |g(E)|^{\frac{1-k^2}{1+k^2}}.
\]
For the map $g$ we use the area distortion inequality
\[ |g(E)| \leq C(K) |E|^{\frac{1}{K}}
\]
from \cite{astalaareadistortion} and conclude the proof by \eqref{eq:1-k}.
\end{proof}

Next, we sketch the reduction of Theorem \ref{thm:riemannmap} to Lemma \ref{lem:areaexpansion}.

\begin{proof}[Proof of Theorem \ref{thm:riemannmap}]
By our assumption $\phi(\infty)=\infty$ and we may also require $\phi(0)=0$.
Since $\phi(\D)$ is a quasidisk $|\phi'(0)| \approx |\phi(1)|$, so we may use the third normalization $\phi(1)=1$.
For an arc $I \subset \partial \D$ let us denote the Carleson square with base $I$ by $Q_I$, that is
$Q_I = \{ z: z/|z| \in I \text{ and } 1-\ell(I) \leq |z| \leq 1 \}$. The  top of $Q_I$ is $z_I=(1-\ell(I))\zeta_I$ where $\zeta_I$ is the center of $I$.
We are going to use the following property of a conformal map to quasidisk target:
for any arc $I$ on $\partial \D$ we have $|\phi'(z_I)| \cdot \diam I \approx \diam \phi(I)$. Furthermore for the top half of the Carleson square $|\phi'(z)| \approx |\phi'(z_I)|$.
We are left to estimate the area of disjoint Carleson squares such that $|\phi'(z_I)| > \rho$. We will only do this in a fixed sector $S$ about $-1$. Let us transfer the situation from the disk to the upper half-plane $\C_+=\{ \Im z >0 \}$. We denote by $T$ the following M\"obius transformation
\[ T(z)=\frac{z-i}{z+i}.
\]
Then $T(\C_+)=\D$ and $T(i)=0$, $T(0)=-1$ and $T(\infty)=1$.
The conjugated map $\psi=T^{-1} \circ \phi \circ T$ is conformal in $\C_+$, globally $K$-quasiconformal
and satisfies $\psi(i)=i$, $\psi(-i)=-i$ and $\psi(\infty)=\infty$.
We choose the inscribed disk $B_j$ inside the image of $Q_I$ under $T^{-1}$ and its reflection along $\R$. Then $\diam(\psi B_j) \gtrsim \rho \diam B_j$ because $|\phi'(z_I)| > \rho$. 
For the set $E = \cup B_j$ we have
\begin{equation*}
 \label{eq:area1}
\rho^2 |E| \approx \sum_{j=1}^n (\rho \diam B_j)^2 \leq \sum_{j=1}^n (\diam \psi B_j)^2 \approx |\psi(E)|.
\end{equation*}
With an appropriate choice of the sector $S$, we may assume $E \subset B(0,1/2)$ and apply 
Lemma \ref{lem:areaexpansion} for the principal map $\psi_1$ with dilatation $\chi_\D \mu_\psi$,
\[ |\psi(E)| \lesssim |\psi_1(E)| \leq C(K) |E|^{1-k}.
\]
Combining the last two inequalities we obtain the desired estimate
\[
 |\{ |\phi'| >\rho \}| \lesssim |\cup_{|\phi'(z_I)| > \rho} Q_I| \approx |E| \lesssim \rho^{-\frac{2}{k}}.
\]
This proves the first part of the Theorem.

In order to prove the second part, we proceed as follows. Consider now an arbitrary Borel set $E \subset \D$, 
\[
 |\phi(E)| = \int_E |\phi'(z)|^2 dm(z) = 2 \int_0^\infty \rho\, |\{z \in E : |\phi'(z)| >\rho \}| d\rho.
\]
We split the integral to two parts at $T=|E|^{-k/2}$. On the interval $[0,T]$ we use the trivial estimate
$T\, |E|$ for the integrand and on $[T,\infty]$ we use the weak-integrability \eqref{eq:weakintegrability}.
The claimed area distortion inequality now follows
\[ |\phi(E)| \leq 2|E| T^2 + C(K) T^{2(1-1/k)} \leq C(K) |E|^{1-k}.
\]
\end{proof}

\bigskip


\bibliographystyle{alpha}
\bibliography{./refer,refextra}

\def\cprime{$'$} \def\cprime{$'$} \def\cprime{$'$}
\begin{thebibliography}{AIM09}

\bibitem[AIM09]{AIM}
Kari Astala, Tadeusz Iwaniec, and Gaven Martin.
\newblock {\em Elliptic partial differential equations and quasiconformal
  mappings in the plane}, volume~48 of {\em Princeton Mathematical Series}.
\newblock Princeton University Press, Princeton, NJ, 2009.

\bibitem[ARS]{astalarohdeschramm95}
Kari Astala, Steffen Rohde, and Oded Schramm.
\newblock Self-similar jordan curves.
\newblock {\em In preparation}.

\bibitem[Ast94]{astalaareadistortion}
Kari Astala.
\newblock Area distortion of quasiconformal mappings.
\newblock {\em Acta Math.}, 173(1):37--60, 1994.

\bibitem[BP88]{becker-pommerenke98}
J.~Becker and Ch. Pommerenke.
\newblock H\"older continuity of conformal maps with quasiconformal extension.
\newblock {\em Complex Variables Theory Appl.}, 10(4):267--272, 1988.

\bibitem[Geh82]{gehringquasidisks}
Frederick~W. Gehring.
\newblock {\em Characteristic properties of quasidisks}, volume~84 of {\em
  S\'eminaire de Math\'ematiques Sup\'erieures}.
\newblock Presses de l'Universit\'e de Montr\'eal, Montreal, Que., 1982.

\bibitem[LSUT]{Lacey-Sawyer-Uriarte}
Michael~T. Lacey, Eric~T. Sawyer, and Ignacio Uriarte-Tuero.
\newblock Astala's conjecture on distortion of {H}ausdorff measures under
  quasiconformal maps in the plane.
\newblock {\em To appear in Acta Math.}

\bibitem[Pom75]{pommerenke75}
Christian Pommerenke.
\newblock {\em Univalent functions}.
\newblock Vandenhoeck \& Ruprecht, G\"ottingen, 1975.
\newblock With a chapter on quadratic differentials by Gerd Jensen, Studia
  Mathematica/Mathematische Lehrb{\"u}cher, Band XXV.

\bibitem[PS09]{prause-smirnov}
Istv\'an Prause and Stanislav Smirnov.
\newblock Quasisymmetric distortion spectrum.
\newblock {\em Preprint, arXiv:0910.4723}, 2009.

\bibitem[Smi]{Smirnov}
Stanislav Smirnov.
\newblock Dimension of quaiscircles.
\newblock {\em To appear in Acta Math.}

\bibitem[Tol09]{tolsaqcdistortion}
Xavier Tolsa.
\newblock Quasiconformal distortion of hausdorff measures.
\newblock {\em Preprint, {arXiv:0907.4933}}, 2009.

\end{thebibliography}

\enlargethispage{2cm}
\end{document}